\documentclass[11pt,letterpaper]{article}

\usepackage[T1]{fontenc}
\usepackage[utf8]{inputenc}
\usepackage{lmodern}
\usepackage[margin=1in]{geometry}
\usepackage{microtype}
\usepackage[hidelinks]{hyperref}
\hypersetup{pdftitle={Lessons from Chess}}
\hypersetup{
  pdftitle={Lessons from Chess: Mathematics in the Age of AI},
  pdfauthor={Dmitriy Leykekhman}
}

\title{Lessons from Chess: Mathematics in the Age of AI}
\author{Dmitriy Leykekhman}
\date{November 23, 2026}

\begin{document}
\maketitle

\section{Introduction}

Rapid advances in artificial intelligence have put mathematics in an unusual position. For generations, mathematics
has been regarded as a distinctly human intellectual activity. 
AI has challenged that view, raising questions about
mathematical research, education, and the future of the profession.
How should we assess students when AI can complete their assignments?
What should constitute an original Ph.D. thesis when powerful AI
assistance is readily available? How will these developments affect
the demand for mathematicians and the organization of mathematical
research?

The mathematical community is divided between anxiety about
the future of the profession and excitement about new possibilities.
Some describe the situation as a crisis; others anticipate a golden
age of mathematics. Universities, funding agencies, and journals
face difficult decisions about how to respond, while the technology
continues to develop faster than its consequences can be understood.
Alongside these technological changes, institutional uncertainty
is already apparent: Dartmouth, for example, is not accepting
applications to its mathematics graduate program for the
2027--28 admissions cycle~\cite{DartmouthMathematics}.

As mathematics confronts this uncertain future, the history
of chess offers a useful point of comparison.
The World Chess Championship match currently taking place provides a timely occasion to reflect on that history.
Chess has already experienced the transition from human supremacy
to overwhelming computer superiority.
Yet human chess has survived and, in many respects, flourished.
What changed, what endured, and what lessons might this experience
offer mathematics?

The chess analogy is not new. It arose, for example, in a recent
conversation between Eric Weinstein and Richard
Borcherds~\cite[at 1:00:25]{BorcherdsWeinstein}.
This opinion note explores the analogy from a historical
perspective, examining both its possibilities and its limitations.
The aim is to draw lessons from the experience of chess
and suggest possible directions for mathematics and its community
as they adapt to the age of AI.

\section{Chess and Mathematics: Parallels}

Chess and mathematics have much in common. Of the first six world
champions, two were mathematicians. The second world champion,
Emanuel Lasker, earned his doctorate in mathematics at the University
of Erlangen in 1900, with a dissertation titled
\emph{\"Uber Reihen auf der Convergenzgrenze}~\cite{Lasker1901}.
The fifth world champion, Max Euwe, earned his doctorate at the
University of Amsterdam in 1926, with a dissertation on differential
invariants~\cite{Euwe1926}.
In 1964, he was appointed professor at Tilburg University.
The sixth world champion, Mikhail Botvinnik, earned a Doctor of
Technical Sciences degree in electrical engineering in 1951.
His research on excitation control and the static stability of
synchronous machines had been published in
1950~\cite{Botvinnik1950}. He was also awarded an honorary degree
in mathematics by the University of Ferrara on September~7, 1991.

As chess became more professional, succeeding in both disciplines
became much harder. Thus, Anatoly Karpov, who originally entered
the mechanics and mathematics program at Moscow State University,
soon had to transfer to the less demanding economics program at
Leningrad University. Surprisingly, Viktor Korchnoi, his future
nemesis, helped him with the transfer~\cite{KarpovGordon2007}.
Of course, there are still individuals accomplished in both
mathematics and chess. For example, John Nunn earned his doctorate
in mathematics at Oxford in 1978 for his work in algebraic
topology~\cite{Nunn1978, Nunn1979} and became British chess champion in 1980.
Both intellectual activities are still highly regarded outside
their respective communities. In the second half of the twentieth
century, both chess and mathematics attracted considerable attention,
surprisingly for much the same reason: politics, particularly
the Cold War.

The popularity of chess in the West reached its peak with the rise
of Bobby Fischer, whose career took off in the late 1950s.
He became the eleventh world champion in 1972 after single-handedly
defeating the Soviet chess machine. The second and third waves of
chess popularity came with the Karpov--Korchnoi and Karpov--Kasparov
battles. These pitted Karpov, a representative of the Soviet system,
against the defector Korchnoi in the first case and the
``child of change'' Kasparov \cite{Kasparov1987} in the second.
The world championship matches between 1972 and 2000 attracted
considerable attention both within and beyond the chess world.

Mathematics also grew in prominence during the Cold War,
as governments poured substantial funding into basic science
and mathematics. In the United States, the National Defense
Education Act of 1958 explicitly supported mathematics education
and provided graduate fellowships and student
loans~\cite{NDEA1958}. The annual number of doctorates awarded
in mathematics and statistics increased from 160 in 1949--50
to approximately 2,000 in 2019--20~\cite{NCESMathDoctorates}.
The expansion was institutional as well as numerical:
the number of U.S. institutions awarding mathematics doctorates
more than doubled, rising from 71 in 1950--1954 to 159
in 1970--1974~\cite[Table~39]{NRC1978}.
Individual departments grew as well. Cornell, my alma mater,
awarded approximately three mathematics doctorates annually
during 1940--1959, compared with more than seven annually
during the 1960s. Of the 73 doctorates awarded in that decade,
52 came in its second half~\cite{CornellGraduateHistory}.
The professional and research infrastructure also expanded:
SIAM was founded in 1952~\cite{SIAMHistory}, while MSRI
in Berkeley and the Institute for Mathematics and its
Applications at Minnesota were established in 1982 \cite{IMAHistory}.

The parallels continue with the rise of computers. The first
chess programs date back to the 1940s. Even the sixth world
champion, Botvinnik, worked on a chess-playing machine and
predicted that computers would become stronger than humans
at chess.
In 1997, Garry Kasparov lost his famous match against Deep Blue
by a score of $3.5$--$2.5$. After the defeat, Kasparov accused
the computer team of cheating, openly claiming that the computer
had not chosen all its moves independently but had occasionally
received human assistance. Sounds familiar?
From a purely chess perspective, the defeat was probably premature;
ten years later, however, humans were no match for computers.
Immediately after the match, everybody started talking about
the death of chess---again. Such discussions go back to the
mid-1920s, when the third world champion, Jos\'e Capablanca,
who rarely lost a game, was already warning about the
``draw death'' of chess~\cite{Capablanca1926}.

Fast forward thirty years, and we find a rather interesting
and surprising situation: chess is more popular than ever.
Millions of people play online every day, including many children.
Chess is no longer a mystery. In the past, only grandmasters---the
chess gurus---could explain a particular move or position,
and even they often made mistakes. Now almost any chess program
can give an accurate evaluation of any position within seconds.
Chess has become much younger and faster. There are fewer
tournaments with classical time controls, but many more rapid
and blitz events. Younger players dominate the field:
both participants in the current world championship match
are twenty years old.
There are numerous video courses, both free and rather pricey,
devoted to particular openings, pawn structures, endgames,
or chess skills. There are also many popular chess streamers,
ranging from the world's elite to beginners and catering
to audiences at every level.
For example, Hikaru Nakamura, one of the world's top players,
currently has 3.7 million followers~\cite{GMHikaru}.
Many platforms broadcast major events live, with commentary
in different languages for audiences around the world.

On the other hand, almost nobody pays attention to the World
Chess Championship match taking place right now. Very few
people outside India, Uzbekistan, and the chess community
even know the players' names.

\section{Chess and Mathematics: Differences}

Of course, chess is not mathematics. Norbert Wiener makes
a rather insightful observation~\cite[p.~21]{Wiener1953}:

\begin{quote}
``The advantage is that mathematics is a field in which one's blunders
tend to show very clearly and can be corrected or erased with a stroke
of the pencil. It is a field which has often been compared with chess,
but differs from the latter in that it is only one's best moments that
count and not one's worst. A single inattention may lose a chess game,
whereas a single successful approach to a problem, among many which have
been relegated to the wastebasket, will make a mathematician's
reputation.''
\end{quote}

Indeed, one can make ten brilliant moves in a row and lose
the game with a single blunder. Magnus Carlsen, the world's
No.~1 player for the past fifteen years, owes his greatness
less to brilliant attacking play than to how rarely he makes
mistakes, combined with his versatility across different
positions and stages of the game.
 In mathematics, on the other hand,
one brilliant idea can make a mathematician famous, despite
less distinguished work elsewhere, even if that work contains mistakes.

In chess, no matter how many games are played or how many opening
variations are investigated, these efforts will produce no theoretical
foundation for quantum mechanics or anything else of value
to other fields.
Chess, however, has the advantage that the value of a new idea
or variation can be tested over the board. In mathematics, the value
of a new idea or development is less clear. It sometimes takes
a long time for new ideas to gain acceptance within the mathematical
community, and some remain completely forgotten.

In addition, mathematics does not have the luxury that chess had:
time to adjust. Progress in computer chess was much slower than
predicted. For example, Botvinnik expected chess programs to play
at master level by 1975; see Chapter~2, ``Chess and Computers,''
of Alexander's \emph{A Book of Chess}~\cite{Alexander1973}.\footnote{%
C.~H.~O'D.~Alexander was an interesting personality in his own right.
He was twice British chess champion and worked at Bletchley Park
during the Second World War on the Enigma project.}
Botvinnik was off by more than a decade.
Only in the late 1980s
did chess computers reach the level of international masters.
After that, their playing strength increased very rapidly.
Mathematics is at a considerable disadvantage compared with chess
in this respect. Despite all their disagreements about the role
of AI in mathematics, mathematicians agree on one thing:
nobody expected such advances in AI even a year ago.

\section{What to Expect Next}

After Deep Blue defeated Garry Kasparov in 1997, it never played
another game. Instead, it was dismantled, and IBM moved on to a new
project. We should therefore expect that, very soon, the high-tech
industry will leave mathematics alone and move on.

There is a nonzero chance that some tech companies will try to
make money from the rapidly growing body of proven mathematical
results. In chess, companies such as ChessBase~\cite{ChessBase27}
sell databases containing essentially every game they can find,
together with efficient tools for searching by position or theme.
There is absolutely no reason not to expect some companies to
do the same for mathematics. All data have value---even mathematical
results---so perhaps this is a good opportunity for arXiv to enter
the game.

However, even if we assume that all the AI and internet tools
we enjoy today will remain available, we still need to address
a question: what new directions should mathematics take?

\section{New Possible Directions for Mathematics}

The mathematical community is not uniform. Mathematics is often
divided into pure, applied, and computational mathematics.
The boundaries are not very clear, but, broadly speaking,
the situation differs considerably across these areas.

Real applied mathematics---mathematics that actually addresses
real-world problems---seems to be in no danger at the moment.
Hard real-world problems are inherently difficult and remain so
even with AI assistance. If AI helps solve one such problem,
great: we can move on to the next one. Models can always be made
more sophisticated and realistic, with more parameters and
additional complexities.
The situation in computational mathematics is similar.
No matter how fast computers become, there will always be
problems that demand their computational power.

Pure (recreational) mathematics seems to be a different story:
this is the area that appears to be in crisis now.
Governments and universities have so far provided substantial
funding for pure mathematics, but this may change very soon.
The common belief that, if even 1\% of pure mathematical research
proves useful, it will pay for everything else is unsupported
and will no longer suffice as a justification.
The question will arise, and very soon: what are pure
mathematicians being paid for? What is the product?
In the past, one could point to a large body of publications,
proved results, books, and citations to justify the funding.
Now, with AI assistance, almost anybody can produce these. 
Human mathematicians are still better at developing theories,
but this advantage may not last much longer either.
So what should we do?

The short answer is to be realistic and adjust.
There are plenty of new opportunities. Unlike chess, mathematics
is enormously broad. Mathematicians with the breadth
of Euler or Poincar\'e have become increasingly rare as the subject
has expanded.
Mathematics now consists of many fields with little interaction.
Even within a single mathematics department, research groups
barely talk to one another or understand one another's work.

One advantage of AI is that it is trained on \textbf{all} available
mathematical literature and can draw on different fields
to establish a particular result.
Some areas of mathematics, such as probability and partial
differential equations, are closely related. For others,
such as algebraic geometry and functional analysis,
the connections are less clear.
AI is capable of finding such connections if they exist.
It would also be interesting to see whether the same result
can be established using tools from different fields.
Historically, results of this kind have been highly prized:
the Atiyah--Singer index theorem connected analysis and topology,
while monstrous moonshine connected finite group theory
and modular functions.
It is a safe bet that many more such connections remain
to be discovered.
If mathematics develops in this direction, the mathematician
of the future should have knowledge spanning a broader range
of mathematics, rather than specialize in a single field---a
return, in some sense, to the mathematicians of the past.

Another idea is inspired by AlphaZero, a successor to AlphaGo.
Starting with the rules of chess and without training on human
games, AlphaZero developed its playing strategy through
reinforcement learning, playing millions of games against itself.
Using substantial computing resources, it reached superhuman
chess strength rather quickly~\cite{Silver2018AlphaZero}.

A natural analogy would be to train an AI to explore a formal
axiom system using specified rules of inference, without supplying
human-written proofs.
For example, an interesting experiment would be to give an AI
Euclid's postulates as a starting point and investigate what
mathematics it develops, including which additional assumptions
it finds necessary.
This is also interesting from a purely philosophical perspective.
Is our human mathematics the only possibility, or could
a completely different, alien mathematics exist?

Education is also an essential part of mathematics. Mathematical
education has traditionally been regarded as valuable, both as
preparation for a wide range of professions and as a way to develop
precise reasoning and problem-solving skills. With AI now able
to perform many tasks that students spend years learning,
this value can no longer be taken for granted.
The mathematical community will have to work hard to attract
and retain mathematics majors and convince the public that
a mathematical education remains worthwhile.
This will require rethinking both what we teach and how we explain
its value: what should students learn, and what will their
mathematical training enable them to do in a world with AI? As the role of research diminishes, the role of education could---and perhaps should---increase.

Finally, we should remember that both chess and mathematics
can be truly beautiful. Fully appreciating that beauty
requires a certain background.
In the past, many magazines and newspapers published chess
puzzles or studies. Even \emph{The New York Times} had a Tuesday
chess column; its chess column was discontinued in 2014.
In mathematics, the situation is less rosy.
Martin Gardner's celebrated \emph{Scientific American} column,
``Mathematical Games,'' ceased appearing regularly
in the early 1980s.

The situation is now changing. Many YouTube channels popularize
mathematics. One of the most prominent is Grant Sanderson's
3Blue1Brown~\cite{3Blue1Brown}, which produces high-quality videos.
Advanced mathematical lectures and courses are also freely
available online---for now---including lectures and video courses by Fields
Medalists~\cite{BorcherdsYouTube}.
However, more needs to be done, probably in a more systematic way.
Mathematics should be popularized at every level, including
elementary school.
Ten years ago, while walking through the halls of my children's
elementary school, I was pleasantly surprised to see that most
kindergartners and first graders had listed mathematics as
their favorite subject.
Apparently, the school had introduced educational computer games
in which students had to solve mathematical problems to advance.
I do not think this enthusiasm for mathematics extended beyond
the first several grades of elementary school.
Platforms such as Brilliant~\cite{Brilliant} are very welcome
and should be promoted more widely by schools and universities.

\section{Conclusion}

AI has created many opportunities and will continue to do so.
In trying to promote their products, AI companies have
unintentionally generated considerable public interest
in mathematics. Many popular podcasts and channels that used
to invite political ``experts'' to comment on world affairs
have started inviting actual mathematicians to explain
what is happening. This is a chance for the mathematical
community to take advantage of the attention and make
the necessary adjustments. Time, however, is short,
and the community needs to act quickly.

The worst thing mathematicians can do now is ignore
the situation or keep declaring that mathematics is dead.
I do not think anyone would enjoy seeing, ten years from now,
mathematics departments at leading universities full of
old people working on problems nobody cares about
and reminiscing about the ``good old days.''

\end{document}